\documentclass[a4paper,10pt]{article}
\usepackage{tikz}
\usepackage{amsfonts,epsf,amsmath,tikz}
\usetikzlibrary{decorations.pathreplacing,automata,calc,positioning}

\newtheorem{theorem}{Theorem}[section]

\newtheorem{proposition}[theorem]{Proposition}

\newtheorem{problem}[theorem]{Problem}

\newcommand{\nex}{$Z$ } 
\newcommand{\qed}{\hfill ~$\square$\bigskip}
\newcommand{\sqed}{\hfill ~$\square$}
\newcommand{\proof}{\noindent{\bf Proof.} }

\newcommand{\ggg}{\gamma_g}

\newcommand{\D}{{Dominator }}
\newcommand{\St}{{Staller }}
\begin{document}

\title{Domination game: effect of edge- and vertex-removal}

\author{
Bo\v{s}tjan Bre\v{s}ar $^a$
\hspace*{-5mm}
\and
Paul Dorbec $^b$
\hspace*{-5mm}
\and
Sandi Klav\v zar $^{c,a}$
\hspace*{-5mm}
\and
Ga\v{s}per Ko\v{s}mrlj $^{c}$
}

\date{}

\maketitle

\begin{center}
$^a$
Faculty of Natural Sciences and Mathematics, University of Maribor, Slovenia\\
{\tt bostjan.bresar@uni-mb.si}
\medskip

$^b$
Univ. Bordeaux, LaBRI, UMR5800, F-33400 Talence \\
CNRS, LaBRI, UMR5800, F-33400 Talence \\
{\tt dorbec@labri.fr}
\medskip

$^c$
Faculty of Mathematics and Physics, University of Ljubljana, Slovenia\\
{\tt sandi.klavzar@fmf.uni-lj.si}\\
{\tt gasper.kosmrlj@student.fmf.uni-lj.si}
\end{center}

\begin{abstract}
The domination game is played on a graph $G$ by two players, named \D and
Staller. They alternatively select vertices of $G$ such that each chosen
vertex enlarges the set of vertices dominated before the move on it. 
Dominator's goal is that the game is finished as soon as possible, 
while \St wants the game to last as long as possible. 
It is assumed that both play optimally. Game 1 and Game 2
are variants of the game in which \D and \St has the first move, respectively.
The game domination number $\gamma_g(G)$, and the Staller-start game domination number
$\gamma_g'(G)$, is the number of vertices chosen in Game 1 and Game 2, respectively.
It is proved that if $e\in E(G)$, then $|\ggg(G) - \ggg(G-e)| \le 2$
and $|\ggg'(G) - \ggg'(G-e)| \le 2$, and that each of the
possibilities here is realizable by connected graphs $G$ for
all values of $\ggg(G)$ and $\ggg'(G)$ larger than 5.
For the remaining small values it is either proved that realizations are not possible or realizing examples are provided.  It is also proved that if $v\in V(G)$, then $\ggg(G) - \ggg(G-v) \le 2$ and $\ggg'(G) - \ggg'(G-v) \le 2$. 
Possibilities here are again realizable by connected graphs $G$ in almost all the cases, the exceptional values are treated similarly as in the edge-removal case.
\end{abstract}


\noindent {\bf AMS Subj. Class. (2010)}: 05C57, 91A43, 05C69

\section{Introduction}

The domination game is played on an arbitrary graph $G$ by two players,
{\em Dominator} and {\em Staller}. They are taking turns choosing
a vertex from $G$ such that whenever they choose a vertex, it dominates
at least one previously undominated vertex. The game ends when all vertices of $G$ are
dominated, so that the set of vertices selected at the end of the game
is a dominating set of $G$.
The aim of \D (Staller) is that the total number of moves played in the
game is as small (as large, resp.) as possible.
By {\em Game 1} ({\em Game 2}) we mean a game in which Dominator (Staller, resp.)
has the first move. Assuming that both players play optimally, the
{\em game domination number} $\gamma_g(G)$ (the {\em
Staller-start game domination number} $\gamma_g'(G)$) of a graph $G$,
denotes the number of vertices chosen in Game 1 (Game 2, resp.).

Note that the domination game is not a combinatorial game in the strict
sense of~\cite{fraenkel-2009}, where the outcome of a game is assumed to
be only of the types (lose, win), (tie, tie) and (draw, draw) for the two players.

The domination game was introduced in~\cite{brklra-2010} (with the idea
going back to~\cite{mike-2003}) and explored by now
from several points of view. Despite the fact that $\gamma(G)\le \ggg(G)\le 2\gamma(G)-1$
holds for any graph $G$ (see~\cite{brklra-2010}), the game domination number is
essentially different from the
domination number. First of all, $\ggg(G)$ is generally much more difficult to
determine than $\gamma(G)$. Even on simple graphs such as paths and cycles, the
problem of determining $\ggg$ is non-trivial~\cite{bill-2012b}.

As proved in~\cite{brklra-2010,bill-2012}, the game domination number and the
Staller-start game domination number can differ only by 1:
$|\gamma_g(G) - \gamma_g'(G)| \le 1$.
Call a pair of integers $(k, \ell)$ {\em realizable} if there exists a graph $G$ with
$\ggg(G) = k$ and $\ggg'(G) = \ell$. Some classes of graphs for possible realizable pairs
are given in~\cite{brklra-2010,brklra-2013,za-2011}. For the complete answer that all
pairs that are potentially realizable can be realized (with relatively simple families
of graphs) see~\cite{gk-2012}.

Kinnersley, West, and Zamani~\cite{bill-2012} conjectured that if $G$ is
an isolate-free forest of order $n$ or an isolate-free graph of order $n$,
then $\gamma_g(G)\le 3n/5$. Actually they posed two conjectures, because
while the truth for isolate-free graphs clearly implies the
truth for isolate-free forests, it is not known
whether the converse implication holds. These conjectures are known as the
$3/5$-conjectures. A progress on them was made in~\cite{brklko-2013}
by constructing large families of trees that attain the conjectured
$3/5$-bound and by finding all extremal trees on up to 20 vertices;
in particular, there are exactly ten trees $T$ on 20 vertices with
$\gamma_g(T)=12$.

Clearly, removing an edge from a graph can only increase its domination number,
that is, $\gamma(G-e) \ge \gamma(G)$. (For an
extensive survey on graphs that are domination critical with respect to edge-
and vertex-removal see~\cite{suwo-1998}.) On the other hand, it was proved
in~\cite{brklra-2013} that for any integer $\ell \geq 1$, there exists a graph
$G$ and its spanning tree $T$ such that $\gamma_g(T)\le \gamma_g(G) - \ell$.
In this paper we answer the question how much $\ggg(G)$ and $\ggg'(G)$ can change
if an edge is removed from $G$. The answer is given in Theorem~\ref{thm:main}
which is followed by ten subsections in which each of the possibilities indicated
by the theorem, is shown to be realizable by connected graphs.
We also ask the analogous question for
vertex-removal and present the answer in Theorem~\ref{thm:vertex}. Again, all
possibilities can be realized by connected graphs. We conclude the paper with
some natural open problems, concerning extensions or generalizations of the
results from this paper.

For a vertex subset $S$ of a graph $G$, let $G|S$ denote the graph $G$ in which vertices from $S$ are considered as being already dominated. In particular, if $S=\{x\}$ we write $G|x$. 
For all the other standard notions not defined in this paper see the monograph on graph domination \cite{hahe-1998}.

In the rest of this section we state some known results to be used in the sequel.

\begin{theorem}
{\rm (\cite[Lemma 2.1]{bill-2012} - Continuation Principle)}
\label{thm:continue}
Let $G$ be a graph and $A, B\subseteq V(G)$.
If $B\subseteq A$, then $\gamma_g(G|A) \le \gamma_g(G|B)$ and
$\gamma_g'(G|A) \le \gamma_g'(G|B)$.
\end{theorem}

\begin{theorem}
{\rm (\cite{brklra-2010,bill-2012})}
\label{thm:we-and-bill}
If $G$ is any graph, then $|\gamma_g(G) - \gamma_g'(G)| \le 1$.
\end{theorem}

\begin{theorem}\label{thm:paths-cycles}
{\rm (\cite{bill-2012b})}
If $n\ge 3$, then
\begin{eqnarray*}
\label{eq:pot1}\ggg(C_n)=
\gamma_g(P_n)&=&\left\{
\begin{array}{ll}\lceil \frac{n}{2}\rceil-1;& n\equiv3 \bmod 4,\\\lceil \frac{n}{2}\rceil;& {\rm otherwise}\,. \end{array}\right.\\
\label{eq:pot2}\gamma_g'(P_n)&=&\left\lceil \frac{n}{2}\right\rceil\,.\\
\label{eq:cikel2}\gamma_g'(C_n)&=&\left\{\begin{array}{ll}\lceil \frac{n-1}{2}\rceil-1;
& n\equiv2 \bmod 4,\\\lceil \frac{n-1}{2}\rceil;& {\rm otherwise}\,. \end{array}
\right.
\end{eqnarray*}
\end{theorem}

\begin{theorem} {\rm\cite[Theorem 4.6]{bill-2012}} \label{thm:drva}
Let $F$ be a forest and $S\subseteq V(F)$. Then $\gamma_g(F|S) \le \gamma_g'(F|S)$.
\end{theorem}

\section{Edge removal}

\begin{theorem}
\label{thm:main}
If $G$ is a graph and $e\in E(G)$, then
$$|\ggg(G) - \ggg(G-e)| \le 2 \qquad {\mbox and} \qquad |\ggg'(G) - \ggg'(G-e)| \le 2\,.$$
\end{theorem}

\proof
To prove the bound $\ggg(G-e) \le \ggg(G) + 2$ it suffices to show that \D has
a strategy on $G-e$ such that at most $\ggg(G) + 2$ moves will be played.
His strategy is to play the game on $G-e$ as follows. In parallel to the real
game he is playing an {\em imagined game} on $G$ by copying every move of \St to
this game and responding optimally in $G$. Each response in the imagined game
is then copied back to the real game in $G-e$. Let $e=uv$ and consider the following
possibilities.

Suppose first that neither \St nor \D play on either of $u$ and $v$ in the course of the real game. Then all the moves in both games are legal and so the imagined game on $G$ lasts no more than $\ggg(G)$ moves. (Recall that \D plays optimally on $G$ but \St might not play optimally.) Since the game on $G-e$ uses the same number of moves, we conclude that in this case the number of moves played in the real game is at most $\ggg(G)$.

Assume now that at some point of the game, the strategy of \D on $G$ is to play
a vertex incident with $e$, say $u$, but this move is not legal in the real game.
This can happen only in the case when $v$ is the only vertex in $N_{G}[u]$ not yet
dominated. In this case \D plays $v$ in the real game and by Theorem~\ref{thm:continue}
(that is, by the Continuation Principle), following
the same strategy \D ensures that the game is finished in no more than
$\ggg(G)$ moves.

Assume next that in the course of the game one of the players played
a vertex incident with $e$, say $u$, and that this is a legal move.
This means that, after this move is copied
into the imagined game on $G$, the vertex $v$ is dominated in this game but may
not yet be dominated in the real game.
If all the moves are legal in the real game (played on $G-e$), then after
at most $\ggg(G)$ moves all vertices except maybe $v$ are dominated. Hence
the real game finishes in no more than $\ggg(G) + 1$ moves.
In the other case \St played a move in which only $v$ was newly dominated,
and this is not a legal move in $G$. Let this move of \St in $G-e$ be the
$k$-th move of the game. Note that
after this move of Staller, the sets of dominated vertices are the
same in both games, denote this set with $D$. Since after the $(k-1)$st move it is Staller's turn
in the imagined game, we derive that
\begin{equation}
\label{eq:aaa}
 (k-1) + \ggg'(G|D)\le \ggg(G)\,.
\end{equation}
(This inequality holds because \St did not necessarily play optimally in the
imagined game.)
Now \D does not copy the move of \St into the imagined
game but simply optimally plays the next moves.
Therefore, since the number of moves left to end each of the games is
$\ggg\left((G-e)|D\right)$, we have:
\begin{eqnarray*}
\ggg(G-e) & \le & k  + \ggg\left((G-e)|D\right) \\
          & = & k  + \ggg(G|D) \\
          & \le & k  + \ggg'(G|D) + 1 \qquad (\mbox{by Theorem~\ref{thm:we-and-bill}})\\
          & \le & \ggg(G) + 2\,. \qquad (\mbox{by \eqref{eq:aaa}})
\end{eqnarray*}

We have thus proved that $\ggg(G-e) \le \ggg(G) + 2$. Note that in the above proof
of this inequality it does not matter whether Game 1 or Game 2 is played on $G-e$.
Hence analogous arguments also give us $\ggg'(G-e) \le \ggg'(G) + 2$.

For the rest of the proof let $A=N_G[u]$.

We next want to demonstrate that  $\ggg(G) \le \ggg(G-e) + 2$. The strategy of
\D on $G$ is to first play on $u$. Then we get
\begin{eqnarray*}
\ggg(G) & \leq & 1 + \ggg'(G|A)\\
		& = & 1 + \ggg'((G-e)|A)\\
		& \leq & 1 + \ggg'(G-e) \qquad (\mbox{by the Continuation Principle})\\
		& \leq & \ggg(G-e) + 2\,. \qquad (\mbox{by Theorem~\ref{thm:we-and-bill}})
\end{eqnarray*}
Note that the equality in the above computation holds because $u$ and $v$ are both in $A$, hence dominated in $G|A$. Thus the edge $e=uv$ is not relevant for any later move.

To complete the proof we show that $\ggg'(G) \le \ggg'(G-e) + 2$. Suppose that \St played first on one of the end vertices of $e$, say $u$. Then we argue as follows:
\begin{eqnarray*}
\ggg'(G) & = & 1 + \ggg(G|A)\\
		& = & 1 + \ggg((G-e)|A)\\
		& \leq & 1 + \ggg(G-e) \qquad (\mbox{by the Continuation Principle})\\
		& \leq & \ggg'(G-e) + 2\,. \qquad (\mbox{by Theorem~\ref{thm:we-and-bill}})
\end{eqnarray*}
Assume now that the first selected vertex $x$ by \St is neither $u$ nor $v$. Then \D replies with the move on $u$. Now we get
\begin{eqnarray*}
\ggg'(G) & = & 1 + \ggg(G|N[x])\\
		& \leq & 2 + \ggg'(G|(N[x]\cup A))\\
		& = & 2 + \ggg'((G-e)|(N[x]\cup A))\\
		& \leq & \ggg'(G-e) + 2\,. \qquad (\mbox{by the Continuation Principle})
\end{eqnarray*}
\qed

In the remainder of this section we demonstrate that all possibilities indicated in Theorem~\ref{thm:main} are
realizable by presenting infinite families of connected graphs for each of the cases.
Two graphs will frequently appear in our constructions, notably $C_6$ and the graph $Z$ from Fig.~\ref{fig:z}.
Recall that $\ggg(C_6)=3=\ggg(C_6|z)$ and $\ggg'(C_6)=2=\ggg'(C_6|z)$, where $z$ is an arbitrary vertex of $C_6$.
Note also that $\ggg(Z)=4=\ggg(Z|z)$ and $\ggg'(Z)=3=\ggg'(Z|z)$.

\begin{figure}[ht!]
\begin{center}
\begin{tikzpicture}[scale=1.0,style=thick]
\def\vr{2pt} 
of vertices
\path (0,0) coordinate (v1); \path (2,0) coordinate (v2);
\path (4,0) coordinate (v3); \path (6,0) coordinate (v4);
\path (5,1) coordinate (v5); \path (0,2) coordinate (v6);
\path (2,2) coordinate (v7); \path (4,2) coordinate (v8);
\path (6,2) coordinate (v9);
\draw (v1) -- (v2) -- (v3) -- (v4) -- (v9) -- (v8) -- (v7) -- (v6) -- (v1);
\draw (v7) -- (v2);
\draw (v8) -- (v5) -- (v4);
\draw (v5) -- (v9);
\draw (v3) -- (v8);
\draw (v1)  [fill=white] circle (\vr); \draw (v2)  [fill=white] circle (\vr);
\draw (v3)  [fill=white] circle (\vr); \draw (v4)  [fill=white] circle (\vr);
\draw (v5)  [fill=white] circle (\vr); \draw (v6)  [fill=white] circle (\vr);
\draw (v7)  [fill=white] circle (\vr); \draw (v8)  [fill=white] circle (\vr);
\draw (v9)  [fill=white] circle (\vr);
\draw [right] (v9) node {$z$};
\end{tikzpicture}
\end{center}
\caption{Graph \nex}
\label{fig:z}
\end{figure}

\subsection{$\ggg(G)-\ggg(G-e)=-2$}\label{sub:edge-2}

\begin{proposition}\label{prop:edge-2}
For any $\ell\geq 3$ there exists a graph $G$ with an edge $e$ such that $\ggg(G)=\ell$ and $\ggg(G-e)=\ell+2$.
\end{proposition}

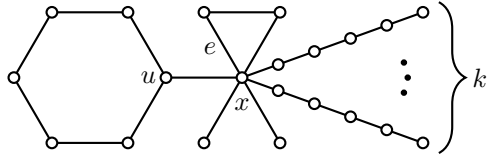
\begin{figure}[ht!]
\begin{center}
\begin{tikzpicture}[scale=1.0,style=thick]
\def\vr{2pt} 
\def\d{1}
\def\x{{2*cos(60)*\d+\d}}
\def\f{20}
\def\dd{2.53209/5}

\path ({cos(60)*\d},{1 + sin(60)*\d}) coordinate (v1);
\path ({cos(60)*\d+\d},{1 + sin(60)*\d}) coordinate (v2);
\path (\x,1) coordinate (v3);
\path ({cos(60)*\d},{1-sin(60)*\d}) coordinate (v4);
\path ({cos(60)*\d+\d},{1-sin(60)*\d}) coordinate (v5);
\path (0,1) coordinate (v6);
\path ({2*cos(60)*\d+2*\d},1) coordinate (v7);
\path ({2*cos(60)*\d+2*\d + cos(120)},{1 + sin(60)*\d}) coordinate (v8);
\path ({3*cos(60)*\d+2*\d},{1 + sin(60)*\d}) coordinate (v9);
\path ({2*cos(60)*\d+2*\d + cos(120)},{1 - sin(60)*\d}) coordinate (v10);
\path ({3*cos(60)*\d+2*\d},{1 - sin(60)*\d}) coordinate (v11);

\path ({2*cos(60)*\d+2*\d + \dd*cos(\f)},{1+\dd*sin(\f)}) coordinate (v12);
\path ({2*cos(60)*\d+2*\d + 2*\dd*cos(\f)},{1+2*\dd*sin(\f)}) coordinate (v13);
\path ({2*cos(60)*\d+2*\d + 3*\dd*cos(\f)},{1+3*\dd*sin(\f)}) coordinate (v14);
\path ({2*cos(60)*\d+2*\d + 4*\dd*cos(\f)},{1+4*\dd*sin(\f)}) coordinate (v15);
\path ({2*cos(60)*\d+2*\d + 5*\dd*cos(\f)},{1+5*\dd*sin(\f)}) coordinate (v16);

\path ({2*cos(60)*\d+2*\d + \dd*cos(\f)},{1-\dd*sin(\f)}) coordinate (v17);
\path ({2*cos(60)*\d+2*\d + 2*\dd*cos(\f)},{1-2*\dd*sin(\f)}) coordinate (v18);
\path ({2*cos(60)*\d+2*\d + 3*\dd*cos(\f)},{1-3*\dd*sin(\f)}) coordinate (v19);
\path ({2*cos(60)*\d+2*\d + 4*\dd*cos(\f)},{1-4*\dd*sin(\f)}) coordinate (v20);
\path ({2*cos(60)*\d+2*\d + 5*\dd*cos(\f)},{1-5*\dd*sin(\f)}) coordinate (v21);
\draw (v1) -- (v2) -- (v3) -- (v5) -- (v4) -- (v6) -- (v1);
\draw (v3) -- (v7) -- (v10);
\draw (v7) -- (v9) -- (v8);
\draw (v8) -- (v7);
\draw (v7) -- (v11);
\draw (v7) -- (v12) -- (v13) -- (v14) -- (v15) -- (v16);
\draw (v7) -- (v17) -- (v18) -- (v19) -- (v20) -- (v21);
\draw (v1)  [fill=white] circle (\vr); \draw (v2)  [fill=white] circle (\vr);
\draw (v3)  [fill=white] circle (\vr); \draw (v4)  [fill=white] circle (\vr);
\draw (v5)  [fill=white] circle (\vr); \draw (v6)  [fill=white] circle (\vr);
\draw (v7)  [fill=white] circle (\vr); \draw (v8)  [fill=white] circle (\vr);
\draw (v9)  [fill=white] circle (\vr);
\draw (v10)  [fill=white] circle (\vr); \draw (v11)  [fill=white] circle
(\vr); \draw (v12)  [fill=white] circle (\vr);
\draw (v13) [fill=white] circle (\vr);
\draw (v14)  [fill=white] circle (\vr);
\draw (v15)  [fill=white] circle (\vr); \draw (v16)
[fill=white] circle (\vr);
\draw (v17)  [fill=white] circle (\vr);
\draw (v18)  [fill=white] circle (\vr); \draw (v19)
[fill=white] circle (\vr);
\draw (v20)  [fill=white] circle (\vr); \draw (v21)
[fill=white] circle (\vr);
\draw [left] (v3) node {$u$};
\draw (v7) node [yshift=-10pt] {$x$};
\draw ({2*cos(60)*\d+2*\d + cos(120)+.08},{1 + sin(60)*\d/2-.05}) node {$e$};
\draw ({2*cos(60)*\d+2*\d + 5*\dd*cos(\f)+.35-.55},1)  [fill=black] circle (0.8pt);
\draw ({2*cos(60)*\d+2*\d + 5*\dd*cos(\f)+.35-.6},1.2)  [fill=black] circle (0.8pt);
\draw ({2*cos(60)*\d+2*\d + 5*\dd*cos(\f)+.35-.6},.8)  [fill=black] circle (0.8pt);
\draw ({2*cos(60)*\d+2*\d + 5*\dd*cos(\f)+.75}, 1) node {$k$};
\draw [decorate,decoration={brace,amplitude=10pt, mirror},xshift=-4pt,yshift=0pt]({2*cos(60)*\d+2*\d + 5*\dd*cos(\f)+.35},0) -- ({2*cos(60)*\d+2*\d + 5*\dd*cos(\f)+.35},2);
\end{tikzpicture}
\end{center}
\caption{Graphs $U_k$}
\label{fig:u_k}
\end{figure}

\proof
We present two different infinite families $U_k$ and $V_k$ realizing odd and even $\ell$, respectively. Let $B$ be the graph isomorphic to $K_{1,4}$ plus an edge, and denote its central vertex by $x$. Let $U_0$ be the graph obtained from the disjoint union of $C_6$ and $B$ by connecting an arbitrary vertex $u$ of the $6$-cycle to $x$ in $B$. The graph $U_k$, $k\geq 1$, is obtained from $U_0$ by identifying one end vertex of $k$ copies of $P_6$ with $x$, see Fig.~\ref{fig:u_k}.

We claim that $\ggg(U_k)=2k+3$ and $\ggg(U_k-e)=2k+5$, where $e$ is one of the two edges incident to $x$ for which $U_k-e$ remains connected, see Fig.~\ref{fig:u_k} again. By Theorem~\ref{thm:main} it suffices to prove that $\ggg(U_k)\leq 2k+3$ and $\ggg(U_k-e)\geq 2k+5$.

For the first inequality we present a strategy for \D that guarantees at most $2k+3$ moves are played on $U_k$. \D starts by playing $x$. Then, he follows \St in the $6$-cycle and in each of the $k$ attached paths, ensuring two moves in each of the subgraphs. Thus $\ggg(U_k)\leq 2k+3$.

To prove the second inequality, Staller's strategy is, if possible, not to be the first to play in the 6-cycle. Note that at least 2 moves will be played in each of the $k$ attached paths, and together with additional 2 moves that will be played in the subgraph that corresponds to $B-e$, it sums up to $2k+2$ moves. If exactly $2k+2$ moves are played before a move in the 6-cycle is played, then it is Dominator who plays first in the 6-cycle, yielding 3 additional moves. Otherwise, at least $2k+3$ moves were played elsewhere, and two additional moves will be played in the 6-cycle, and so $\ggg(U_k-e)\geq 2k+5$. This concludes the proof for the case when $\ell$ is odd.
\medskip

For the case when $\ell$ is even we construct a family $V_k$ in a similar way as $U_k$ by replacing the $6$-cycle with the graph \nex from Fig.~\ref{fig:z}. More precisely for $V_0$ we take the disjoint union of $Z$ and $B$, and add an edge connecting $z$ and $x$. Then the graph $V_k$, $k\geq 1$, is obtained from $V_0$ by identifying one end vertex of $k$ copies of $P_6$ with $x$. By using parallel arguments to the above case one can derive that $\ggg(V_k)=2k+4$ and $\ggg(V_k-e)=2k+6$.
\qed

To round off the subsection we show that when $\ell <3$, there exists no graph $G$ such that
$\ggg(G) = \ell$ and $\ggg(G-e)=\ell-2$. Indeed, we have:
\begin{itemize}
\item If $\ggg(G) = 1$, then $\ggg(G-e)\le 2$.\\
If $\ggg(G)=1$, then $G$ has a universal vertex $v$. In $G-e$, \D can play
$v$ and in this way dominate all but at most one vertex $w$. Hence \St has to dominate $w$ in any legal move. \sqed
\item If $\ggg(G)=2$, then $\ggg(G-e)\le 3$.\\
We prove this by contradiction: suppose $\ggg(G)=2$ but $\ggg(G-e)=4$. We propose a strategy
of \St in $G$ that will require at least 3 moves to be played. Suppose \D plays on some vertex $d_1$ in $G$, and let $s_1,d_2,s_2$ be the next three optimal moves played in the game on $G-e$. Let $s'_1, d'_2, s'_2$ be vertices (not necessarily distinct from $s_1, d_2, s_2$) that are newly dominated in $G-e$ when $s_1,d_2,s_2$ are played, respectively. Note that in $G$, if \St plays on $s_1$, at most one of $d'_2$ and $s'_2$ may be dominated. Thus since $\ggg(G)=2$, the move $s_1$ is not a legal answer. Therefore the edge $e$ necessarily connects $s'_1$ either to $d_1$ or to $s_1$. Now if \St plays on $d_2$, she does not dominate $s'_2$, hence $\ggg(G)\ge 3$, a contradiction. \sqed
\end{itemize}

\subsection{$\ggg(G)-\ggg(G-e)=-1$}

Let $T$ be a tree with $\ggg(T)=\ell$, $\ell\geq1$, and let $v$ be an optimal start vertex for Dominator. Let $G_\ell$ be the graph obtained from $T$ by attaching two additional leaves to $v$ and identifying $v$ with a vertex of a triangle. Note that $|V(G_\ell)|=|V(T)|+4$. Let $e$ be an edge of the triangle incident with $v$, having $y$ as the other end vertex. By the Continuation Principle, \D will not play on any of the new four vertices added to the tree $T$. Hence $v$ is also an optimal start vertex for \D in $G_\ell$, and so $\ggg(G_\ell)=\ell$.

In $G_\ell-e$, \D starts by playing $v$, and will play $y$ only if it is the only legal move. If \St plays on $y$ or its neighbor at the $k$th move, then \D continues optimally and the total number of moves is $\ggg(T|D) + k$, where $D$ is the set of dominated vertices at this stage in $T$. Applying Theorem~\ref{thm:drva}, this is at most $\ggg'(T|D) + k = \ggg(T)+1=\ell +1$. Otherwise, dominator plays $y$ on his last move and the game also ends after $\ell+1$ moves.
To see that $\ggg(G_\ell-e)\ge \ell+1$ we present Staller's strategy. Whenever \D plays in the subgraph of $G_\ell-e$ that corresponds to $T$, she responds in this subgraph as well, by playing as if the game was played in $T$.  Note that if \D plays the neighbor of $y$, then \St will be the first to play in the remainder of the game with respect to $T$, which will thus in total take at least $\ell$ steps, by using Theorem \ref{thm:drva}. Hence at least $\ell+1$ moves will be played in $G_\ell-e$. If \D does not play on the neighbor of $y$ during the game, then in the last move \St plays on it, which concludes the proof of $\ggg(G_\ell-e)=\ell+1$. \sqed

\subsection{$\ggg(G)-\ggg(G-e)=0$}

From Theorem~\ref{thm:paths-cycles} we get that $\ggg(C_n)=\ggg(C_n-e)$ for any $n\geq 3$. \sqed

\subsection{$\ggg(G)-\ggg(G-e)=1$}\label{sub:edge1}

\begin{proposition}\label{prop:edge+1}
For any $\ell\geq 3$ there exists a graph $G$ with an edge $e$ such that $\ggg(G)=\ell$ and $\ggg(G-e)=\ell-1$.
\end{proposition}

\proof
For the case when $\ell=3$, take the disjoint union of two complete graphs of order at least three and add two edges that form a matching in the resulting graph. Clearly, the game domination number of this graph is $3$, while after removing one of the two additional edges the game domination number drops to $2$.

\begin{figure}[ht!]
\begin{center}
\begin{tikzpicture}[scale=1.0,style=thick]
\def\vr{2pt} 
of vertices
\path (0,0) coordinate (v1); \path (0,2) coordinate (v2);
\path (1,1) coordinate (v3); \path (1,2) coordinate (v4);
\path (1,3) coordinate (v5); \path (2,0) coordinate (v6);
\path (2,2) coordinate (v7);
\path ({2+cos(60)},{2+sin(60)}) coordinate (v8);
\path ({2-cos(60)},{2+sin(60)}) coordinate (v9);
\path (3,2) coordinate (v10);
\path (4,2) coordinate (v11);
\path ({2+cos(-30)},{2+sin(-30)}) coordinate (v12);
\path ({2+2*cos(-30)},{2+2*sin(-30)}) coordinate (v13);
\draw (v1) -- (v2) -- (v4) -- (v7)-- (v6) -- (v1);
\draw (v3) -- (v7) -- (v8);
\draw (v7) -- (v9);
\draw (v4) -- (v5);
\draw (v11) -- (v10) -- (v7) -- (v12) -- (v13);
\draw (v2) -- (v3);
\draw (v1)  [fill=white] circle (\vr); \draw (v2)  [fill=white] circle (\vr);
\draw (v3)  [fill=white] circle (\vr); \draw (v4)  [fill=white] circle (\vr);
\draw (v5)  [fill=white] circle (\vr); \draw (v6)  [fill=white] circle (\vr);
\draw (v7)  [fill=white] circle (\vr); \draw (v8)  [fill=white] circle (\vr);
\draw (v9)  [fill=white] circle (\vr);
\draw (v10)  [fill=white] circle (\vr); \draw (v11)  [fill=white] circle
(\vr); \draw (v12)  [fill=white] circle (\vr);
\draw (v13) [fill=white] circle (\vr);
\draw (v7) node [yshift=-8pt, xshift=4pt] {$x$};
\draw (v2) node [above] {$x'$};
\draw [below] (v3) node {$y$};
\draw (v4) node [below] {$t$};
\draw (0.35,1.4) node {$e$};
\draw (3.49,1.79)  [fill=black] circle (0.8pt);
\draw (3.5,1.6)  [fill=black] circle (0.8pt);
\draw (3.44,1.43)  [fill=black] circle (0.8pt);
\draw (4.5, 1.35) node {$k$};
\draw [decorate,decoration={brace,amplitude=10pt, mirror},xshift=-4pt,yshift=0pt]({4.3-.3},0.8) -- ({4.3},2.1);
\end{tikzpicture}
\end{center}
\caption{Graphs $Y_k$}
\label{fig:y_k}
\end{figure}
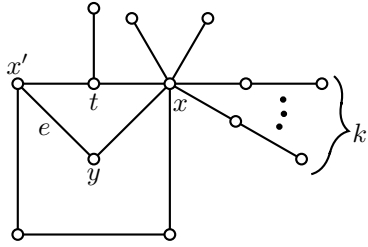

For the general case when $\ell\geq4$, we construct the following family of graphs denoted by $Y_k$, $k\geq0$. Let $Y_0$ be obtained in the following way. Let $t$ be a vertex of some $C_5$ and $x,x'$ its neighbors. Add a new vertex $y$ connected to both $x$ and $x'$. Finally, attach two leaves to $x$ and one to $t$. For $k\geq 1$, the graph $Y_k$ is obtained from $Y_0$ by identifying the end vertices of $k$ copies of $P_3$ with $x$, see Fig.~\ref{fig:y_k}. We claim that $\ggg(Y_k)=k+4$ and $\ggg(Y_k-e)=k+3$.

Note that if \D plays his first move on $x$, then only $k+3$ vertices remain undominated which already yields $\ggg(Y_k)\leq k+4$. Next we present the strategy for \St which ensures that at least $k+4$ moves are needed to end the game in $Y_k$.

If \D starts on $x$, then \St responds on $y$. Then there are still $k+2$ isolated vertices left undominated, no pair of which has a common neighbor in $Y_k$. Hence $k+4$ moves will be played in total. Otherwise, if \D does not start on $x$, then \St responds on a leaf adjacent to $x$. It follows that in the subgraph of $Y_k$ that corresponds to $Y_0$ at least 4 moves will be played. In turn at least $k$ moves will be played in the $k$ attached paths, thus at least $k+4$ moves are needed. We conclude that $\ggg(Y_k)=k+4$.

To prove that $\ggg(Y_k-e)\leq k+3$ we first explain the strategy of Dominator. He starts on $x$. If \St dominates two vertices in the next move, then $k+1$ vertices remain undominated and the bound is ensured. Otherwise, in his second move \D dominates two new vertices by playing either $t$ or $x'$. This gives the desired upper bound for $\ggg(Y_k-e)$. Finally, we present a strategy for \St that guarantees at least $k+3$ moves will be played in $Y_k-e$. If \D does not play $x$, then the strategy for \St is the same as above when the game was played in $Y_k$ (in particular, she responds on a leaf adjacent to $x$). Otherwise, if \D plays $x$ in his first move, then \St responds on $t$. Then $k+1$ isolated vertices, no pair of each has a common neighbor, remain undominated, yielding $\ggg(Y_k-e)\geq k+3$. This concludes the proof.
\qed

There exists no graph $G$ with $\ggg(G)=2$ and $\ggg(G-e)=1$. Actually, we have the following:
\begin{itemize}
\item If $\ggg(G)\ge 2$, then $\ggg(G-e)\ge 2$. If $\ggg(G)\ge 4$, then $\ggg(G-e)\ge 3$.\\
It is proved in~\cite{brklra-2010} that $\ggg(G)\le 2\gamma(G) - 1$, or, equivalently,
$\gamma(G)\ge \frac{\ggg(G)+1}{2}$. Together with the fact that the domination number does not decrease by edge removal, this implies
\begin{equation}
\ggg(G-e)\ge \gamma(G-e)\ge \gamma(G) \ge \left\lceil \frac{\ggg(G)+1}{2}\right\rceil\label{eq:no-small-ell}
\end{equation}
and this yields the desired bounds for $\ggg(G)=2,3$ or 4. \sqed
\end{itemize}

\subsection{$\ggg(G)-\ggg(G-e)=2$}\label{sub:edge2}

\begin{proposition}\label{prop:edge+2}
For any $\ell\geq 5$ there exists a graph $G$ with an edge $e$ such that $\ggg(G)=\ell$ and $\ggg(G-e)=\ell-2$.
\end{proposition}

\proof
We present families of graphs $X_k$ and $Q_k$ that realize even and odd values of $\ell$, respectively. The arguments for the first family are given in full detail, while the reasoning for $Q_k$ is analogous.

\begin{figure}[ht!]
\begin{center}
\begin{tikzpicture}[scale=0.9,style=thick]
\def\vr{2pt} 
of vertices
\path (0,0) coordinate (v1); \path (2,0) coordinate (v2);
\path (4,0) coordinate (v3); \path (6,0) coordinate (v4);
\path (5,1) coordinate (v5); \path (0,2) coordinate (v6);
\path (2,2) coordinate (v7); \path (4,2) coordinate (v8);
\path (6,2) coordinate (v9);
\path (7,1) coordinate (v10);
\path (8,1) coordinate (v11); \path (8,2) coordinate (v12);
\path ({8+cos(144)},{2+sin(144)}) coordinate (v13);
\path ({8+cos(108)},{2+sin(108)}) coordinate (v14);
\path ({8+cos(72)},{2+sin(72)}) coordinate (v15);
\path ({8+cos(36)},{2+sin(36)}) coordinate (v16);
\path (9,2) coordinate (v17);
\path (10,2) coordinate (v18);
\path (11,2) coordinate (v19);
\path (12,2) coordinate (v20);
\path (13,2) coordinate (v21);

\path ({8+cos(-20)},{2+sin(-20)}) coordinate (v22);
\path ({8+2*cos(-20)},{2+2*sin(-20)}) coordinate (v23);
\path ({8+3*cos(-20)},{2+3*sin(-20)}) coordinate (v24);
\path ({8+4*cos(-20)},{2+4*sin(-20)}) coordinate (v25);
\path ({8+5*cos(-20)},{2+5*sin(-20)}) coordinate (v26);

\draw (v1) -- (v2) -- (v3) -- (v4) -- (v9) -- (v8) -- (v7) -- (v6) -- (v1);
\draw (v7) -- (v2);
\draw (v8) -- (v5) -- (v4);
\draw (v5) -- (v9); \draw (v3) -- (v8);
\draw (v5) -- (v10) -- (v9);
\draw (v4) -- (v10) -- (v8);
\draw (v10) -- (v11);
\draw (v9) -- (v12) -- (v11);
\draw (v13) -- (v12) -- (v14);
\draw (v15) -- (v12) -- (v16);
\draw (v26) -- (v25) -- (v24) -- (v23) -- (v22) -- (v12) -- (v17) -- (v18) -- (v19) -- (v20) -- (v21);

\draw (v1)  [fill=white] circle (\vr); \draw (v2)  [fill=white] circle (\vr);
\draw (v3)  [fill=white] circle (\vr); \draw (v4)  [fill=white] circle (\vr);
\draw (v5)  [fill=white] circle (\vr); \draw (v6)  [fill=white] circle (\vr);
\draw (v7)  [fill=white] circle (\vr); \draw (v8)  [fill=white] circle (\vr);
\draw (v9)  [fill=white] circle (\vr);
\draw (v10)  [fill=white] circle (\vr); \draw (v11)  [fill=white] circle
(\vr); \draw (v12)  [fill=white] circle (\vr);
\draw (v13) [fill=white] circle (\vr);
\draw (v14)  [fill=white] circle (\vr);
\draw (v15)  [fill=white] circle (\vr); \draw (v16)
[fill=white] circle (\vr);
\draw (v17)  [fill=white] circle (\vr);
\draw (v18)  [fill=white] circle (\vr); \draw (v19)
[fill=white] circle (\vr);
\draw (v20)  [fill=white] circle (\vr); \draw (v21)
[fill=white] circle (\vr);
\draw (v22)  [fill=white] circle (\vr);
\draw (v23) [fill=white] circle (\vr);
\draw (v24) [fill=white] circle (\vr);
\draw (v25) [fill=white] circle (\vr);
\draw (v26) [fill=white] circle (\vr);
\draw [above] (v9) node {$z$};
\draw [above] (v10) node[xshift=3pt] {$z'$};
\draw (v12) node[xshift=-5pt, yshift=-5pt] {$x$};
\draw (7.5,0.8)node {$e$};
\draw[right] (v11) node {$x'$};
\draw (12.3,1.47)  [fill=black] circle (0.8pt);
\draw (12.3,1.25)  [fill=black] circle (0.8pt);
\draw (12.24,1.05)  [fill=black] circle (0.8pt);
\draw (13.7, .95) node {$k$};
\draw [decorate,decoration={brace,amplitude=10pt, mirror},xshift=-4pt,yshift=0pt]({8+5*cos(-20)+.3},{2+5*sin(-20)-.3}) -- (13.4,2.1);
\end{tikzpicture}
\end{center}
\caption{Graphs $X_k$}
\label{fig:x_k}
\end{figure}
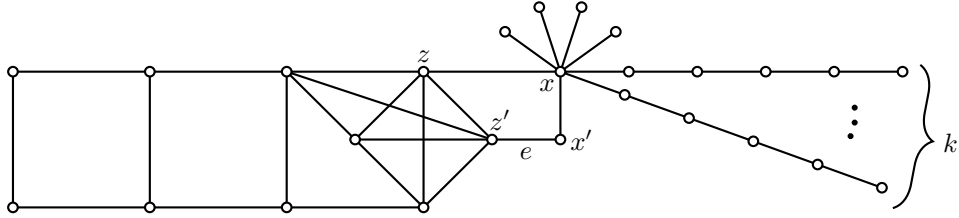

We construct $X_0$ as follows. Duplicate the vertex $z$ in $Z$ (see Fig.~\ref{fig:z}) obtaining a new vertex $z'$ with the same closed neighborhood as $z$, and denote the resulting graph by $Z'$. Next, take the disjoint union of $Z'$ with $K_{1,5}$ having $x$ as its center and denote one of its leaves by $x'$. Finally we get $X_0$ by connecting $z$ with $x$, and $z'$ with $x'$. The graph $X_k$, $k\geq 1$, is obtained from $X_0$ by identifying one end vertex of $k$ copies of $P_6$ with $x$, see Fig.~\ref{fig:x_k}. We set $e$ to be the edge between $z'$ and $x'$.

We claim that $\ggg(X_k)=2k+6$ and $\ggg(X_k-e)=2k+4$. By Theorem~\ref{thm:we-and-bill} it suffices to present a strategy for \D yielding $\ggg(X_k-e)\leq 2k+4$, and a strategy for \St which gives $\ggg(X_k)\geq 2k+6$. To show the first inequality, \D starts the game by playing $x$. Any move of \St in one of the $k$ attached paths is followed by a move of \D in the same path, so that all vertices of this path are dominated. With this strategy \St is forced to be the first to play in the subgraph that corresponds to $Z'$. Since $\ggg'(Z')=3$, \D can ensure that only three moves are played in this subgraph. Altogether we get that $2k+4$ moves will be played in $X_k-e$.

It remains to present a strategy for \St in $X_k$. Whenever \D plays on one of the $k$ attached paths, \St follows on the same path in such a way that all vertices on the path at distance at least 2 from $x$ are dominated after her move. If \D plays one of the vertices $x$ or $x'$, \St responds with a move on the other vertex from $\{x,x'\}$, if this is possible. Note that this is not possible only in the case when $z'$ was dominated before. By this strategy, \St forces \D to be the first to play in the subgraph isomorphic to $Z'$. Suppose first that when \D starts to play in $Z'$, $x$ and $x'$ have already been played. Since $\ggg(Z')=4$, four moves will be played in $Z'$, hence together with $2k$ moves on the attached paths the total number of moves sums up to $2k+6$. Otherwise, if $x$ and $x'$ have not been played at the time when \D starts to play in $Z'$, then \St responds by playing on one of the leaves attached to $x$. If the next move of Dominator which is not played on one of the attached paths, is also played in $Z'$, then \St responds by playing on one the leaves attached to $x$, again. Since $\gamma(Z')=3$, at this point in the game there are still undominated vertices left in $Z'$ as well as two undominated leaves attached to $x$. Thus at least two more moves are needed, altogether at least $2k+6$ moves. On the other hand, if the next move of Dominator which is not played on one of the attached paths, is played on $x$, then Staller's next move is in $Z'$ ensuring four moves will be played in $Z'$. Again we get that at least $2k+6$ moves will be played in $X_k$ in total which concludes the proof for even $\ell$.

\begin{figure}[ht!]
\begin{center}
\begin{tikzpicture}[scale=0.9,style=thick]
\def\vr{2pt} 
of vertices
\path (6.5,2) coordinate (v9);
\path (6.5,1) coordinate (v10);
\path (v9) +(90+360/7:1) coordinate (v4);
\path (v4) +(90+360*2/7:1) coordinate (v5);
\path (v5) +(90+360*3/7:1) coordinate (v6);
\path (v6) +(90+360*4/7:1) coordinate (v7);
\path (v7) +(90+360*5/7:1) coordinate (v8);
\path (8,1) coordinate (v11); \path (8,2) coordinate (v12);
\path ({8+cos(144)},{2+sin(144)}) coordinate (v13);
\path ({8+cos(108)},{2+sin(108)}) coordinate (v14);
\path ({8+cos(72)},{2+sin(72)}) coordinate (v15);
\path ({8+cos(36)},{2+sin(36)}) coordinate (v16);
\path (9,2) coordinate (v17);
\path (10,2) coordinate (v18);
\path (11,2) coordinate (v19);
\path (12,2) coordinate (v20);
\path (13,2) coordinate (v21);

\path ({8+cos(-20)},{2+sin(-20)}) coordinate (v22);
\path ({8+2*cos(-20)},{2+2*sin(-20)}) coordinate (v23);
\path ({8+3*cos(-20)},{2+3*sin(-20)}) coordinate (v24);
\path ({8+4*cos(-20)},{2+4*sin(-20)}) coordinate (v25);
\path ({8+5*cos(-20)},{2+5*sin(-20)}) coordinate (v26);

\draw (v4) -- (v5) -- (v6) -- (v7) -- (v8) -- (v9) -- (v10) -- (v4);
\draw (v9) -- (v4);
\draw (v8) -- (v10);
\draw (v10) -- (v11);
\draw (v9) -- (v12) -- (v11);
\draw (v13) -- (v12) -- (v14);
\draw (v15) -- (v12) -- (v16);
\draw (v26) -- (v25) -- (v24) -- (v23) -- (v22) -- (v12) -- (v17) -- (v18) -- (v19) -- (v20) -- (v21);

\draw (v4)  [fill=white] circle (\vr);
\draw (v5)  [fill=white] circle (\vr); \draw (v6)  [fill=white] circle (\vr);
\draw (v7)  [fill=white] circle (\vr); \draw (v8)  [fill=white] circle (\vr);
\draw (v9)  [fill=white] circle (\vr);
\draw (v10)  [fill=white] circle (\vr); \draw (v11)  [fill=white] circle
(\vr); \draw (v12)  [fill=white] circle (\vr);
\draw (v13) [fill=white] circle (\vr);
\draw (v14)  [fill=white] circle (\vr);
\draw (v15)  [fill=white] circle (\vr); \draw (v16)
[fill=white] circle (\vr);
\draw (v17)  [fill=white] circle (\vr);
\draw (v18)  [fill=white] circle (\vr); \draw (v19)
[fill=white] circle (\vr);
\draw (v20)  [fill=white] circle (\vr); \draw (v21)
[fill=white] circle (\vr);
\draw (v22)  [fill=white] circle (\vr);
\draw (v23) [fill=white] circle (\vr);
\draw (v24) [fill=white] circle (\vr);
\draw (v25) [fill=white] circle (\vr);
\draw (v26) [fill=white] circle (\vr);
\draw [above] (v9) node {$z$};
\draw [below] (v10) node {$z'$};
\draw (v12) node[xshift=-5pt, yshift=-5pt] {$x$};
\draw (7.25,1.2)node {$e$};
\draw[right] (v11) node {$x'$};
\draw (12.3,1.47)  [fill=black] circle (0.8pt);
\draw (12.3,1.25)  [fill=black] circle (0.8pt);
\draw (12.24,1.05)  [fill=black] circle (0.8pt);
\draw (13.7, .95) node {$k$};
\draw [decorate,decoration={brace,amplitude=10pt, mirror},xshift=-4pt,yshift=0pt]({8+5*cos(-20)+.3},{2+5*sin(-20)-.3}) -- (13.4,2.1);
\end{tikzpicture}
\end{center}
\caption{Graphs $Q_k$}
\label{fig:q_k}
\end{figure}
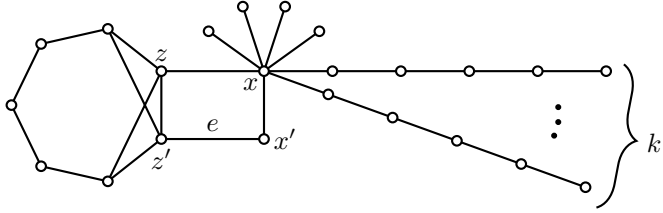

The family $Q_k$ which realizes the case when $\ell$ is odd is constructed as follows. Take a copy of $C_6$, denote one of its vertices by $z$, and add a duplicate vertex $z'$ of $z$ so that this two vertices have the same closed neighborhoods. Denote the resulting graph by $Z''$ and take the disjoint union of $Z''$ with $K_{1,5}$ having $x$ as its center, and denote one of its leaves by $x'$. Finally we get $Q_0$ by connecting $z$ with $x$, and $z'$ with $x'$. The graph $Q_k$, $k\geq 1$, is obtained from $Q_0$ by attaching $k$ copies of $P_6$ at their end vertices to $x$, see Fig.~\ref{fig:q_k}. We set $e$ to be the edge between $z'$ and $x'$. As noted in the beginning of the proof, the arguments for $\ggg(Q_k)=2k+5$ and $\ggg(Q_k-e)=2k+3$ follow similar lines as above.
\qed

By inequality~\eqref{eq:no-small-ell}, there exists no graph $G$ such that $\ggg(G)=\ell$ and $\ggg(G-e)=\ell-2$
for some edge $e$ when $\ell\le 4$.

\subsection{$\ggg'(G)-\ggg'(G-e)=-2$}

Similarly as in Subsection~\ref{sub:edge-2}, one can verify that $\ggg'(U_k)=2k+4$ and $\ggg'(U_k-e)=2k+6$ for any $k\geq0$. Also, $\ggg'(V_k)=2k+5$ and $\ggg'(V_k-e)=2k+7$. In particular, note that the optimal first move of \St is to play on a leaf adjacent to $x$. Hence:

\begin{proposition}\label{prop:edge-2-prime}
For any $\ell\geq 4$ there exists a graph $G$ with an edge $e$ such that $\ggg'(G)=\ell$ and $\ggg'(G-e)=\ell+2$.
\end{proposition}

Note also that for $\ell<4$, there are no graphs such that $\ggg'(G)=\ell$ and $\ggg'(G-e)=\ell+2$ for some edge $e$ in $G$. Indeed, we have:
\begin{itemize}
\item If $\ggg'(G) = 1$, then $\ggg'(G-e)=2$.\\
Obviously, the only non-trivial graphs $G$ with $\ggg'(G)=1$ are complete graphs, and $\ggg'(K_n-e)=2$. \sqed
\item If $\ggg'(G)=2$, then $\ggg'(G-e)\le 3$.\\
Suppose $\ggg'(G)=2$. For any move $s_1$ of \St in $G$, \D has an answer $d_1$ that dominates all of $G$. 
This move played in $G-e$ would thus dominate all but at most one vertex of $G-e$.
The next move of \St has to dominate that vertex and thus $\ggg'(G-e)\le 3$. \sqed
\item If $\ggg'(G)=3$, then $\ggg'(G-e)\le 4$.\\
Suppose $\ggg'(G)=3$. Let $s_1$ be an optimal move of \St in $G-e$, and let \D answer to $s_1$
by the same move as in $G$. Consider an optimal reply $s_2$ of \St in $G-e$. If $s_2$ is legal in $G$, then $\{s_1,d_1,s_2\}$ is a dominating set of $G$, so it dominates all $G-e$ but at most one vertex, and hence $\ggg'(G-e)\le 4$. If $s_2$ is not legal in $G$, it means that $s_2$ newly dominates only one end of $e$, the other end being $s_1$ or $d_1$. After Staller's move $s_2$ in $G-e$, the set of dominated vertices is then exactly the same as the set of dominated vertices after the two first moves in $G$ and both ends of $e$ are dominated, so any legal move in $G-e$ finishes the game, and again $\ggg'(G-e)\le 4$. \sqed
\end{itemize}

\subsection{$\ggg'(G)-\ggg'(G-e)=-1$}

By Theorem~\ref{thm:paths-cycles}, $\ggg'(C_{2\ell+1})=\ell$ while $\ggg'(C_{2\ell+1}-e)=\ell+1$ for any $\ell\geq 1$. \sqed

\subsection{$\ggg'(G)-\ggg'(G-e)=0$}

Note first that $\ggg'(C_4)=\ggg'(C_4-e)=2$. Let $G$ be the graph obtained from $P_4$ by identifying one of its inner vertices denoted by $u$ with a vertex of a triangle. Then $\ggg'(G)=\ggg'(G-e)=3$ where $e$ is the edge of the triangle not incident with $u$. Let $k\geq 4$ and let $G_k$ be the graph obtained from $K_k$ by attaching one leaf to every vertex of $K_k$. Let $e$ be an edge of $G_k$ that lies in the $k$-clique. Then it is straightforward that $\ggg'(G_k)=k=\ggg'(G_k-e)$. \sqed

\subsection{$\ggg'(G)-\ggg'(G-e)=1$}

Similarly as in Subsection~\ref{sub:edge1}, one can verify that $\ggg'(Y_k)=k+5$ and $\ggg'(Y_k-e)=k+4$ for any $k\geq0$. In particular, note that the optimal first move of \St is to play on a leaf adjacent to $x$. Consider next the graph $H$ obtained from the disjoint union of $K_{1,4}$ and a triangle by joining with an edge the center of $K_{1,4}$ with one vertex of the triangle, and by adding edge $e$ between another vertex of the triangle and a leaf of $K_{1,4}$. Then $\ggg'(H)=4$ and $\ggg'(H)=3$. Hence we have:

\begin{proposition}\label{prop:edge1-prime}
For any $\ell\geq 4$ there exists a graph $G$ with an edge $e$ such that $\ggg'(G)=\ell$ and $\ggg'(G-e)=\ell-1$.
\end{proposition}

Note that when $\ell < 4$, there are no graphs with $\ggg'(G)=\ell$ and $\ggg'(G-e)=\ell -1$ for some edge $e$. Indeed:
\begin{itemize}
\item If $G$ is a graph with at least one edge, then $\ggg'(G-e)\ge 2$.\\
This is clear because $G-e$ is not complete. \sqed
\item If $\ggg'(G)=3$, then $\ggg'(G-e)\ge 3$.\\
Let $s_1$ be an optimal first move of \St in $G$. There are no vertices in $G$ that dominate
all of $V(G)\setminus N[s_1]$, hence there are none either in $G-e$. Thus the game in
$G-e$ requires at least two more moves and $\ggg'(G-e)\ge 3$. \sqed
\end{itemize}

\subsection{$\ggg'(G)-\ggg'(G-e)=2$}

Similarly as in Subsection~\ref{sub:edge2}, one can verify that $\ggg'(X_k)=2k+7$ and $\ggg'(X_k-e)=2k+5$ for any $k\geq0$. Also, $\ggg'(Q_k)=2k+6$ and $\ggg'(Q_k-e)=2k+4$. In particular, note that the optimal first move of \St is to play on a leaf adjacent to $x$. Hence:

\begin{proposition}
For any $\ell\geq 6$ there exists a graph $G$ with an edge $e$ such that $\ggg'(G)=\ell$ and $\ggg'(G-e)=\ell-2$.
\end{proposition}

When $\ell < 6$, there are no graphs with $\ggg'(G)=\ell$ and $\ggg'(G-e)=\ell -2$ for some edge $e$. Indeed,
\begin{itemize}
\item for $\ell < 4$ we proved the assertion in the previous subsection.
\item If $\ggg'(G)=4$, then $\ggg'(G-e)\ge 3$.\\
Suppose by way of contradiction that $\ggg'(G-e) = 2$. Then to any move $s_1$ of \St in $G$, \D answers $d_1$ as if in $G-e$. Since $\{s_1,d_1\}$ dominates $G-e$, it also dominates $G$. \sqed
\item If $\ggg'(G)=5$, then $\ggg'(G-e)\ge 4$.\\
Suppose by way of contradiction that $\ggg'(G-e) = 3$. Then to any move $s_1$ of \St in $G$, \D answers $d_1$ as if in $G-e$. Now any legal move $s_2$ of \St in $G-e$ dominates $G-e$ and so dominates $G$. However, \St may play a move in $G$ that was not legal in $G-e$, but then \D answering with $s_2$ finishes the game in at most 4 moves. \sqed
\end{itemize}

\section{Vertex removal}

In contrast to the fact that $\ggg(G-e) \le \ggg(G) + 2$ holds, the game domination number of a vertex-deleted subgraph of $G$ cannot be bounded above by a function of the game domination number of $G$. This is not surprising because the same phenomenon holds for the usual domination number (and because $\gamma(G)\leq\ggg(G)\leq 2\gamma(G)-1$). More explicitly, let $k$ be a non-negative integer and let $H$ be an arbitrary graph with $\ggg(H)=k+1$. Let $G$ be the graph obtained from $H$ by adding to it a universal vertex $v$. Then $\ggg(G)=1$ and hence $\ggg(G-v)-\ggg(G)=k$. The same construction works for the Staller-start game domination number.

On the other hand, we prove the following:

\begin{theorem}\label{thm:vertex}
If $G$ is a graph and $v\in V(G)$, then
$$\ggg(G) - \ggg(G-v) \le 2 \qquad {\mbox and} \qquad \ggg'(G) - \ggg'(G-v) \le 2\,.$$
\end{theorem}

\proof
To prove the first inequality, let \D start on $v$ when Game 1 is played in $G$.
We get
\begin{eqnarray*}
\ggg(G) & \leq & 1 + \ggg'(G|N[v])\\
		& = & 1 + \ggg'((G-v)|(N[v]-\{v\}))\\
		& \leq & 1 + \ggg'(G-v) \qquad (\mbox{by the Continuation Principle})\\
		& \leq & \ggg(G-v) + 2\,. \qquad (\mbox{by Theorem~\ref{thm:we-and-bill}})
\end{eqnarray*}
In Game 2 we consider two cases. In the first case \St plays $v$ in her first move. We get that
\begin{eqnarray*}
\ggg'(G) & = & 1 + \ggg(G|N[v])\\
		& = & 1 + \ggg((G-v)|(N[v]-\{v\}))\\
		& \leq & 1 + \ggg(G-v) \qquad (\mbox{by the Continuation Principle})\\
		& \leq & \ggg'(G-v) + 2\,. \qquad (\mbox{by Theorem~\ref{thm:we-and-bill}})
\end{eqnarray*}
In the second case \St chooses a vertex $x$, $x\neq v$. Then \D responds by playing $v$, hence
\begin{eqnarray*}
\ggg'(G) & = & 1 + \ggg(G|N[x])\\
		& \leq & 2 + \ggg'(G|(N[x]\cup N[v]))\\
		& = & 2 + \ggg'((G-v)|(N[x]\cup N[v]-\{v\}))\\
		& \leq & \ggg'(G-v) + 2\,. \qquad (\mbox{by the Continuation Principle})
\end{eqnarray*}
\qed

We have already observed, that $\ggg(G-v)-\ggg(G)$ as well as
$\ggg'(G-v)-\ggg'(G)$
can be arbitrarily small. In the rest of this section we construct
infinite families
of (connected) graphs demonstrating that for any $t\in \{0,1,2\}$ and any
integer
$\ell \ge 5$ (or smaller---depending of the case), there exists a graph
$G$ with $\ggg(G) = \ell$ (resp. $\ggg'(G) = \ell$) and $\ggg(G-v)-\ggg(G)
= t$
(resp. $\ggg'(G-v)-\ggg'(G) = t$).

\subsection{$\ggg(G)-\ggg(G-v)=0$}\label{sub:vertex0}

Let $\ell$ be a positive integer and let $G'$ be an arbitrary graph with
$\ggg(G') = \ell$. Let $x$ be an optimal start vertex for Dominator and
let $G$ be the graph obtained from $G'$ by attaching a leaf $v$ to $x$
(actually, we
could attach any number of leaves). We claim that $\ggg(G) = \ggg(G-v)=\ell$.
Clearly, $\ggg(G-v)=\ell$ since $G-v=G'$. By the Continuation principle,
\D would
start the game rather on $v$ than on $x$. But then $x$ is an optimal start
vertex
for \D also in $G$, hence $\ggg(G)=\ell$.
\sqed

\subsection{$\ggg(G)-\ggg(G-v)=1$}\label{sub:vertex1}

To see that for any integer $\ell\ge 2$ there exists a graph $G$ such that $\ggg(G)=\ell$ and $\ggg(G-v)=\ell-1$ for some $v\in V(G)$ it suffices to notice that the sequence $\big(\ggg(P_n)\big)_{n\geq1}$ is unbounded, non-decreasing and $\ggg(P_{n+1})-\ggg(P_n)\leq 1$ for any $n$. \sqed

\subsection{$\ggg(G)-\ggg(G-v)=2$}\label{sub:vertex2}

\begin{proposition}\label{prop:vertex-2}
For any $\ell\geq 5$ there exists a graph $G$ with a vertex $v$ such that $\ggg(G)=\ell$ and $\ggg(G-v)=\ell-2$.
\end{proposition}

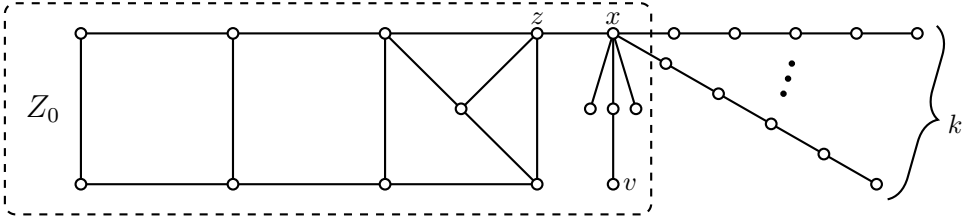
\begin{figure}[ht!]
\begin{center}
\begin{tikzpicture}[scale=1.0,style=thick]
\def\vr{2pt} 
of vertices
\path (0,0) coordinate (v1); \path (2,0) coordinate (v2);
\path (4,0) coordinate (v3); \path (6,0) coordinate (v4);
\path (5,1) coordinate (v5); \path (0,2) coordinate (v6);
\path (2,2) coordinate (v7); \path (4,2) coordinate (v8);
\path (6,2) coordinate (v9);
\path (7,2) coordinate (v10);
\path (7,1) coordinate (v11); \path (7,0) coordinate (v12);
\path (7.8,2) coordinate (v13); \path (8.6,2) coordinate (v14);
\path (9.4,2) coordinate (v15); \path (10.2,2) coordinate (v16); \path (11,2)
coordinate (v17);
\path ({7+cos(-30)*0.8},{2+sin(-30)*0.8}) coordinate (v18);
\path ({7+cos(-30)*1.6},{2+sin(-30)*1.6}) coordinate (v19);
\path ({7+cos(-30)*2.4},{2+sin(-30)*2.4}) coordinate (v20);
\path ({7+cos(-30)*3.2},{2+sin(-30)*3.2}) coordinate (v21);
\path ({7+cos(-30)*4},{2+sin(-30)*4}) coordinate (v22);
\path (6.7,1.0) coordinate (v23);
\path (7.3,1.0) coordinate (v24);
\draw (v1) -- (v2) -- (v3) -- (v4) -- (v9) -- (v8) -- (v7) -- (v6) -- (v1);
\draw (v7) -- (v2);
\draw (v8) -- (v5) -- (v4);
\draw (v5) -- (v9); \draw (v3) -- (v8);
\draw (v9) -- (v10) -- (v11) -- (v12);
\draw (v10) -- (v13) -- (v14) -- (v15) -- (v16) -- (v17);
\draw (v10) -- (v18) -- (v19) -- (v20) -- (v21) -- (v22);
\draw (v24) -- (v10) -- (v23);
\draw (v1)  [fill=white] circle (\vr); \draw (v2)  [fill=white] circle (\vr);
\draw (v3)  [fill=white] circle (\vr); \draw (v4)  [fill=white] circle (\vr);
\draw (v5)  [fill=white] circle (\vr); \draw (v6)  [fill=white] circle (\vr);
\draw (v7)  [fill=white] circle (\vr); \draw (v8)  [fill=white] circle (\vr);
\draw (v9)  [fill=white] circle (\vr);
\draw (v10)  [fill=white] circle (\vr); \draw (v11)  [fill=white] circle
(\vr); \draw (v12)  [fill=white] circle (\vr);
\draw (v13) [fill=white] circle (\vr);
\draw (v14)  [fill=white] circle (\vr);
\draw (v15)  [fill=white] circle (\vr); \draw (v16)
[fill=white] circle (\vr);
\draw (v17)  [fill=white] circle (\vr);
\draw (v18)  [fill=white] circle (\vr); \draw (v19)
[fill=white] circle (\vr);
\draw (v20)  [fill=white] circle (\vr); \draw (v21)
[fill=white] circle (\vr);
\draw (v22)  [fill=white] circle (\vr);
\draw (v23) [fill=white] circle (\vr);
\draw (v24) [fill=white] circle (\vr);
\draw [above] (v9) node {$z$};
\draw [above] (v10) node {$x$};
\draw [right] (v12) node {$v$};
\draw (9.35,1.6)  [fill=black] circle (0.8pt);
\draw (9.30,1.4)  [fill=black] circle (0.8pt);
\draw (9.24,1.2)  [fill=black] circle (0.8pt);
\draw (-.5,1) node {\large $Z_0$};
\draw (11.5, .8) node {$k$};
\draw[rounded corners, dashed](-1,-.4)rectangle(7.5,2.4);
\draw [decorate,decoration={brace,amplitude=10pt, mirror},xshift=-4pt,yshift=0pt](10.75,-.2) -- (11.4,2.1);
\end{tikzpicture}
\end{center}
\caption{Graphs $Z_k$}
\label{fig:z_k}
\end{figure}

\proof
We start the proof for even $\ell$ by presenting the following family $Z_k$, $k\geq 0$.  Let $S$ be the graph obtained from $K_{1,3}$ with $x$ as its center in which one edge is subdivided. Denote the vertex that is not in $N_S[x]$ by $v$. Then $Z_0$ is obtained from the disjoint union of $Z$ and $S$ by connecting $x$ and $z$ with an edge. See Fig.~\ref{fig:z_k} where $Z_0$ is encircled by a dashed curve. The graph $Z_k$, $k\geq 1$, is obtained from $Z_0$ by identifying the end vertex of $k$ copies of $P_6$ with $x$, see Fig.~\ref{fig:z_k} again.

We claim that $\ggg(Z_k)=2k+6$ and $\ggg(Z_k-v)=2k+4$. By Theorem~\ref{thm:vertex} it suffices to show that $\ggg(Z_k)\geq 2k+6$ and $\ggg(Z_k-v)\leq 2k+4$.

To prove the first assertion consider the following strategy of Staller. We first observe that on each of the $k$ attached paths at least two vertices different from $x$ will be played, hence at least $2k$ vertices in total. Moreover, at least two vertices will be played in the subgraph of $Z_k$ that corresponds to $S$. If exactly $2k+2$ moves are played on this part, \St can force \D to be the first one to play in the subgraph isomorphic to $Z$. In this case four moves will be played in this subgraph and hence at least $2k+6$ moves in total. Otherwise, if \St is forced to play first in $Z$, then at least $2k+3$ moves were played on the rest of $Z_k$. Since at least three moves will be played in $Z$ (note that $\gamma(Z)=3$) again at least $2k+6$ moves will be played on $Z_k$.

To prove that  $\ggg(Z_k-v)\leq 2k+4$ consider the strategy of \D to play first on $x$. By following \St in $Z$ and in each of the $k$ attached paths he ensures that three moves will be played in $Z$ and two in each of the $k$ attached paths. Hence in total $2k+4$ moves will be played. This proves the proposition in the case when $\ell$ is even.

We use a similar construction to prove the result for odd $\ell$. In the construction of $Z_k$ replace $Z$ by $C_6$, denoting any of its vertices by $z$. Let the resulting graph be denoted by $W_k$. We claim that
$\ggg(W_k)=2k+5$ and $\ggg(W_k-v)=2k+3$. Note that $\ggg(C_6)=3=\ggg(C_6|z)$ and $\ggg'(C_6)=2=\ggg'(C_6|z)$. Then we argue that this is indeed the case with arguments parallel to those that we used for the graphs $Z_k$.
\qed

Note that there does not exist a graph $G$ such that $\ggg(G)=4$ and $\ggg(G-v)=2$ for a vertex $v$. Indeed after the first optimal move of \D in $G-v$, the set $C$ of undominated vertices induces a complete subgraph of $G-v$, and any vertex in $G-C$ that is adjacent to a vertex of $C$ dominates the entire $C$. In $G$, \D can start by playing on the same vertex so that only vertices of $C\cup \{v\}$ are left undominated. Clearly at most two more moves will be played in $G$, hence $\ggg(G) \leq 3$. It is also easy to see that there does not exist a graph $H$ with a vertex $v$ such that $\ggg(H)=3$ and $\ggg(H-v)=1$.

\subsection{$\ggg'(G)-\ggg'(G-v)=0$}\label{sub:vertex0-prime}

Let $\ell$ be a positive integer and let $G$ be the graph obtained from
$K_{\ell+2}$
by attaching a leaf to $\ell$ of its vertices. $G$ is thus of order
$2\ell+2$.
Let $v$ be one of the two vertices of the clique that has no leaf attached.
Then it is not difficult to see that $\ggg'(G) = \ggg'(G-v) = \ell +1$.
\sqed

\subsection{$\ggg'(G)-\ggg'(G-v)=1$}

One can use paths in the same way as in Subsection~\ref{sub:vertex1}.

\subsection{$\ggg'(G)-\ggg'(G-v)=2$}

Similarly as in Subsection~\ref{sub:vertex2}, one can verify that $\ggg'(Z_k)=2k+7$ and $\ggg'(Z_k-v)=2k+5$ for any $k\geq0$. Also, $\ggg'(W_k)=2k+6$ and $\ggg'(W_k-v)=2k+4$. In particular, note that the optimal first move of \St is to play on a leaf adjacent to $x$.

The graph $H$ obtained from attaching a leaf $v$ to any vertex of $C_6$ provides $\ggg'(H)=4$ and $\ggg'(H-v)=2$. Similarly, for the graph $H'$ obtained from $Z$ by attaching a leaf $v$ to the vertex $z$ we get that $\ggg'(H')=5$ and $\ggg'(H'-v)=3$. Hence we have the following.

\begin{proposition}\label{prop:vertex-2-prime}
For any $\ell\geq 4$ there exists a graph $G$ with a vertex $v$ such that $\ggg'(G)=\ell$ and $\ggg'(G-v)=\ell-2$.
\end{proposition}

\section{Concluding remarks}
We conclude the paper by two problems that arise from the results of this paper.

\begin{problem}
\label{problem1}
Given a positive integer $k$, can one find a general upper and lower bound for $\ggg(G)-\ggg(G_{k})$ where $G_{k}$ is obtained from a graph $G$ by deletion of $k$ edges from $G$?
\end{problem}

An interesting instance of Problem~\ref{problem1} is the question whether
$|\ggg(G)-\ggg(G-\{e,e'\})|$ can be 3 or 4.

\begin{problem}
Which of the subsets of $\{-2,-1,0,1,2\}$ can be realized as $$\{\ggg(G)-\ggg(G-e):\ e\in E(G)\}$$ within the family of all (respectively connected) graphs $G$?
\end{problem}

In particular, does there exist a graph $G$ with edges denoted by $e_{-2},e_{-1},e_{0},e_{1},e_{2}$ such that $\ggg(G)-\ggg(G-e_i)=i$ for all $i$?

In addition, one can ask for a characterization of certain subfamilies of graphs with respect to the above properties. For instance, following domination terminology a possible question is to characterize the graphs that are game domination edge-critical. That is, for which $G$ we have $\{\ggg(G)-\ggg(G-e):\ e\in E(G)\}\subseteq\{-2,-1\}$?

\section*{Acknowledgments}
This work was done in the frame of the bilateral France-Slovenian project BI-FR/13-14-PROTEUS-003 entitled Graph domination.
B.B. and S.K. are supported by the Ministry of Science of Slovenia under the grant P1-0297, and are also with the Institute
of Mathematics, Physics and Mechanics, Jadranska 19, 1000 Ljubljana.
G.K. is also financed part by the European Union - European Social Fund,
and by Ministry of Economic Development and Technology of Republic of
Slovenia.


\begin{thebibliography}{99}

\bibitem{brklra-2010}
  B.~Bre{\v{s}}ar, S.~Klav{\v{z}}ar, D.~F.~Rall,
  Domination game and an imagination strategy,
  SIAM J. Discrete Math. 24 (2010) 979--991.

\bibitem{brklra-2013}
  B.~Bre{\v{s}}ar, S.~Klav{\v{z}}ar, D.~F.~Rall,
  Domination game played on trees and spanning subgraphs,
  Discrete Math. 313 (2013) 915--923.

\bibitem{brklko-2013}
  B.~Bre{\v{s}}ar, S.~Klav{\v{z}}ar, G.~Ko{\v{s}}mrlj, D.~F.~Rall,
  Domination game: extremal families of graphs for the 3/5-conjectures,
  Discrete Appl. Math. 161 (2013) 1308--1316.


\bibitem{fraenkel-2009}
  A.~S.~Fraenkel,
  Combinatorial games: selected bibliography with a succinct gourmet introduction,
  Electron. J. Combin. (August 9, 2012) DS2, 109pp.

\bibitem{hahe-1998}
  T.~W.~Haynes, S.~T.~Hedetniemi, P.~J.~Slater,
  Fundamentals of Domination in Graphs,
  Marcel Dekker, New York, 1998.

\bibitem{mike-2003}
  M.~Henning, personal communication, 2003.

\bibitem{bill-2012}
  W.~B.~Kinnersley, D.~B.~West, R.~Zamani,
  Extremal problems for game domination number,
  manuscript, 2012.

\bibitem{bill-2012b}
  W.~B.~Kinnersley, D.~B.~West, R.~Zamani,
  Game domination for grid-like graphs,
  manuscript, 2012.

\bibitem{gk-2012}
  G.~Ko\v{s}mrlj,
  Realizations of the game domination number,
  to appear in J. Comb. Optim., DOI: 10.1007/s10878-012-9572-x.

\bibitem{suwo-1998}
  D.~P.~Sumner, E.~Wojcicka,
  Graph critical with respect to the domination number,
  Chapter 16 in: Domination in Graphs: Advanced Topics (T.~W.~Haynes, S.~T.~Hedetniemi, P.~J.~Slater, eds.),
    Marcel Dekker, New York, 1998, 439--469.

\bibitem{za-2011}
  R.~Zamani, Hamiltonian cycles through specified edges in bipartite graphs,
  domination game, and the game of revolutionaries and spies,
  Ph. D. Thesis, University of Illinois at Urbana-Champaign.
  Pro-Quest/UMI, Ann Arbor (Publication No. AAT 3496787)

\end{thebibliography}
\end{document}